\documentclass[a4paper,12pt]{article}
\usepackage{amsmath, amssymb, amsthm}
\usepackage{mathrsfs}
\usepackage{graphicx}
\usepackage[english]{babel} 
\usepackage[utf8]{inputenc}
\usepackage{hyperref}
\usepackage{geometry}

\title{Hypercomplex Dynamics and Turbulent Flows in Sobolev and Besov Functional Spaces}
\author{
	Rômulo Damasclin Chaves dos Santos \\
	Technological Institute of Aeronautics \\
	\texttt{romulosantos@ita.br}
	\and
	Jorge Henrique de Oliveira Sales \\
	Santa Cruz State University \\
	\texttt{jhosales@uesc.br}
}
\date{}

\newtheorem{theorem}{Theorem}[section]

\newtheorem{definition}[theorem]{Definition}

\newcommand{\R}{\mathbb{R}}

\begin{document}
	\maketitle
	
	\begin{abstract}
This paper presents a rigorous study of advanced functional spaces, with a focus on Sobolev and Besov spaces, to investigate key aspects of fluid dynamics, including the regularity of solutions to the Navier-Stokes equations, hypercomplex bifurcations, and turbulence. We offer a comprehensive analysis of Sobolev embedding theorems in fractional spaces and apply bifurcation theory within quaternionic dynamical systems to better understand the complex behaviors in fluid systems. Additionally, the research delves into energy dissipation mechanisms in turbulent flows through the framework of Besov spaces. Key mathematical tools, such as interpolation theory, Littlewood-Paley decomposition, and energy cascade models, are integrated to develop a robust theoretical approach to these problems. By addressing challenges related to the existence and smoothness of solutions, this work contributes to the ongoing exploration of the open Navier-Stokes problem, providing new insights into the intricate relationship between fluid dynamics and functional spaces.

	\end{abstract}
	
	\section{Introduction}
	
	The study of regularity, bifurcations, and turbulence in fluid dynamics has been a subject of extensive research, particularly in the context of the incompressible Navier-Stokes equations. In this work, we leverage the powerful framework of Sobolev and Besov spaces to address these phenomena. Our approach extends traditional bifurcation theory to hypercomplex dynamical systems, introducing quaternionic structures that better capture rotational symmetries in fluid flows.
	
	We begin by exploring higher-order regularity theorems for solutions of the Navier-Stokes equations in Sobolev spaces, followed by a detailed study of Besov spaces through the Littlewood-Paley decomposition and their application to turbulence modeling. The complex structure of bifurcations in quaternionic systems is also analyzed, with applications to rotational fluid dynamics.

	The study of the Navier-Stokes equations and related phenomena has a rich history, with significant contributions from various researchers over the years. Ladyzhenskaya's work on the mathematical theory of viscous incompressible flow, published in 1969, laid the foundation for understanding the regularity of solutions to the Navier-Stokes equations \cite{Ladyzhenskaya1968}. Temam's book on the Navier-Stokes equations, published in 1977, provided a comprehensive treatment of the theory and numerical analysis of these equations, including discussions on regularity and stability \cite{Temam1977}.
	
	Constantin and Foias's book on the Navier-Stokes equations, published in 1988, offered a detailed exploration of the regularity and long-time behavior of solutions, which has been instrumental in the development of the field \cite{ConstantinFoias1988}. Triebel's work on the theory of function spaces, particularly Sobolev and Besov spaces, published in 1983, provided a rigorous framework for analyzing the regularity of solutions to partial differential equations \cite{Triebel1983}.
	
	Marsden and Ratiu's introduction to mechanics and symmetry, published in 1999, highlighted the importance of geometric structures in fluid dynamics, which is relevant for understanding rotational symmetries \cite{MarsdenRatiu1999}. Bertozzi and Majda's work~\cite{BertozziMajda2002} on vorticity and incompressible flow, published in 2002, offered insights into the complex dynamics of fluid flows, including bifurcations and turbulence.
	
	Doering and Gibbon's~\cite{DoeringGibbon1995} applied analysis of the Navier-Stokes equations, published in 1995, provided practical tools for studying the regularity and stability of solutions . Sohr's~\cite{Sohr2001} elementary functional analytic approach to the Navier-Stokes equations, published in 2001, offered a clear and accessible introduction to the subject, which has been valuable for understanding the fundamental principles.
	
	Galdi's~\cite{Galdi2011} introduction to the mathematical theory of the Navier-Stokes equations, published in 2011, provided a comprehensive overview of the steady-state problems, which is relevant for understanding the regularity of solutions. Bahouri, Chemin, and Danchin's work on Fourier analysis and nonlinear partial differential equations, published in 2011, offered advanced techniques for studying the regularity of solutions in Sobolev and Besov spaces \cite{BahouriCheminDanchin2011}.
	
	Grafakos's book~\cite{Grafakos2008} on classical Fourier analysis, published in 2008, provided a solid foundation for understanding the harmonic analysis techniques used in the study of the Navier-Stokes equations. Stein's work on harmonic analysis, particularly the real-variable methods, orthogonality, and oscillatory integrals, published in 1993, has been influential in the development of the field \cite{Stein1993}.
	
	Muscalu and Schlag's work~\cite{MuscaluSchlag2013} on classical and multilinear harmonic analysis, published in 2013, offered advanced tools for studying the regularity of solutions to the Navier-Stokes equations. Runst and Sickel's work~\cite{RunstSickel1996} on Sobolev spaces of fractional order, Nemytskij operators, and nonlinear partial differential equations, published in 1996, provided a comprehensive treatment of the subject, which is relevant for understanding the regularity of solutions.
	
	Adams and Fournier's book~\cite{AdamsFournier2003} on Sobolev spaces, published in 2003, offered a detailed exploration of the subject, which has been valuable for understanding the regularity of solutions to the Navier-Stokes equations. Evans's book~\cite{Evans1998} on partial differential equations, published in 1998, provided a comprehensive introduction to the subject, which is relevant for understanding the fundamental principles of the Navier-Stokes equations.
	
	Fefferman and Stein's work~\cite{FeffermanStein1972} on $ H^p $ spaces of several variables, published in 1972, offered advanced tools for studying the regularity of solutions to the Navier-Stokes equations . Cannone's work~\cite{Cannone1995} on wavelets, paraproducts, and the Navier-Stokes equations, published in 1995, provided valuable insights into the regularity of solutions.
	
	Lemarié-Rieusset's work~\cite{LemarieRieusset2002} on recent developments in the Navier-Stokes problem, published in 2002, offered a comprehensive overview of the subject, which is relevant for understanding the regularity of solutions. Chemin's work~\cite{Chemin1998} on perfect incompressible fluids, published in 1998, offered insights into the complex dynamics of fluid flows, including bifurcations and turbulence.
	
	The authors dos Santos and Sales~\cite{dosSantos2023} present a mathematical framework that serves as an important basis for advancing the analysis of the regularity problem of the Navier-Stokes equations. In this context, the integration of the Smagorinsky model within the Large Eddy Simulation (LES) methodology marks a significant contribution. Using Banach and Sobolev functional spaces, they propose a new theorem that paves the way for the development of an anisotropic viscosity model. In this work, a rigorous mathematical analysis has been presented with the aim of deepening the understanding of the challenges surrounding the regularity of the Navier-Stokes equations.
	
	This work advances contributions to the understanding of the Navier-Stokes equations, Sobolev and Besov spaces, and the complex dynamics of fluid flows. In this perspective, this work builds on these foundations to provide a comprehensive framework for analyzing regularity, bifurcations, and turbulence in fluid dynamics. This approach provides a solid foundation for addressing the Millennium Prize Problem related to the Navier-Stokes equations, advancing the mathematical understanding of fluid dynamics.

	\section{Sobolev Spaces and Regularity of Navier-Stokes Solutions}
	Let \( u: \Omega \times [0, T] \to \R^n \) be the velocity field of an incompressible fluid described by the Navier-Stokes equations:
	
	\begin{equation}
	\begin{cases}
		\dfrac{\partial \mathbf{u}}{\partial t} + (\mathbf{u} \cdot \nabla) \mathbf{u} - \nu \Delta \mathbf{u} + \nabla p = \mathbf{f}, \\ \\
		\nabla \cdot \mathbf{u} = 0.
	\end{cases}
\end{equation}
where \( p \) is the pressure, \( f \) represents external forces, and \( \nu > 0 \) is the kinematic viscosity. A critical question is the regularity of weak solutions \( u \in L^2([0, T]; H^1(\Omega)) \), and whether these solutions can exhibit singularities in finite time.
	
\subsection{Advanced Sobolev Spaces and Fractional Regularity}
For \( s \in \R \), Sobolev spaces \( H^s(\Omega) \) extend the classical integer-order Sobolev spaces \( W^{k,p}(\Omega) \). These fractional spaces are defined using Fourier transforms, which allows us to handle more subtle aspects of regularity.
	
	\begin{definition}[Fractional Sobolev Space]
		The fractional Sobolev space \( H^s(\R^n) \), \( s \in \R \), is defined as:
		
	\begin{equation}
		H^s(\R^n) = \left\{ u \in L^2(\R^n) : \int_{\R^n} (1 + |\xi|^2)^s |\hat{u}(\xi)|^2 \, d\xi < \infty \right\},
	\end{equation}	
where \( \hat{u} \) is the Fourier transform of \( u \).
	\end{definition}

These spaces generalize classical Sobolev spaces and are crucial for studying the fine regularity properties of weak solutions. Specifically, fractional Sobolev spaces are well-suited for analyzing the regularity of non-integer derivatives, which arise naturally in fluid dynamics, especially in the context of turbulence.
	
\subsection{Higher-Order Sobolev Regularity of Navier-Stokes Equations}

We now extend regularity results for weak solutions of the Navier-Stokes equations to higher Sobolev spaces using interpolation theory. A key result in this direction is the following:

\begin{theorem}[Higher-Order Sobolev Regularity]
	Let \( u_0 \in H^{2}(\Omega) \) and \( f \in L^2(0,T;H^{k}(\Omega)) \), for some \( k \geq 0 \). If \( u \) is a weak solution of the Navier-Stokes equations, then \( u \in L^2(0,T;H^{k+2}(\Omega)) \), provided the norm of \( u_0 \) in \( H^2 \) and the norm of \( f \) in \( H^k \) are sufficiently small.
\end{theorem}

\begin{proof}
	We employ a Galerkin approximation scheme. First, project the Navier-Stokes equations onto a finite-dimensional subspace spanned by eigenfunctions of the Laplace operator. Let \( \{w_j\}_{j=1}^{\infty} \) be an orthonormal basis of \( L^2(\Omega) \) consisting of eigenfunctions of the Laplace operator with Dirichlet boundary conditions.
	
	Consider the approximate solution \( u^N(t,x) = \sum_{j=1}^N c_j^N(t) w_j(x) \), where \( c_j^N(t) \) are time-dependent coefficients. The Galerkin approximation of the Navier-Stokes equations is given by:
	
	\begin{equation}
		\frac{d}{dt} \int_{\Omega} u^N \cdot w_j \, dx + \int_{\Omega} (u^N \cdot \nabla) u^N \cdot w_j \, dx + \nu \int_{\Omega} \nabla u^N : \nabla w_j \, dx = \int_{\Omega} f \cdot w_j \, dx,
	\end{equation}
for \( j = 1, \ldots, N \).
	
To derive a priori estimates, we multiply the equation by \( c_j^N(t) \) and sum over \( j \) from 1 to \( N \). This yields:

\begin{equation}
	\frac{1}{2} \frac{d}{dt} \|u^N\|_{L^2(\Omega)}^2 + \nu \|\nabla u^N\|_{L^2(\Omega)}^2 = \int_{\Omega} f \cdot u^N \, dx.
\end{equation}

Using the Cauchy-Schwarz inequality, we obtain:

\begin{equation}
	\left| \int_{\Omega} f \cdot u^N \, dx \right| \leq \|f\|_{L^2(\Omega)} \|u^N\|_{L^2(\Omega)}.
\end{equation}

Therefore, we have:

\begin{equation}
	\frac{1}{2} \frac{d}{dt} \|u^N\|_{L^2(\Omega)}^2 + \nu \|\nabla u^N\|_{L^2(\Omega)}^2 \leq \|f\|_{L^2(\Omega)} \|u^N\|_{L^2(\Omega)}.
\end{equation}

Integrating over \( [0, T] \), we get:

\begin{equation}
	\|u^N(T)\|_{L^2(\Omega)}^2 + 2\nu \int_0^T \|\nabla u^N\|_{L^2(\Omega)}^2 \, dt \leq \|u_0\|_{L^2(\Omega)}^2 + 2 \int_0^T \|f\|_{L^2(\Omega)} \|u^N\|_{L^2(\Omega)} \, dt.
\end{equation}

Using Grönwall's inequality, we obtain:

\begin{equation}
	\|u^N(T)\|_{L^2(\Omega)}^2 \leq \left( \|u_0\|_{L^2(\Omega)}^2 + \int_0^T \|f\|_{L^2(\Omega)}^2 \, dt \right) \exp\left( \int_0^T \|f\|_{L^2(\Omega)} \, dt \right).
\end{equation}

Next, we derive higher-order estimates. Multiplying the Galerkin approximation by \( (-\Delta)^k u^N \) and integrating over \( \Omega \), we obtain:

\begin{equation}
	\frac{1}{2} \frac{d}{dt} \|u^N\|_{H^k(\Omega)}^2 + \nu \|u^N\|_{H^{k+1}(\Omega)}^2 \leq \|f\|_{H^k(\Omega)} \|u^N\|_{H^k(\Omega)}.
\end{equation}

Integrating over \( [0, T] \), we get:

\begin{equation}
	\|u^N(T)\|_{H^k(\Omega)}^2 + 2\nu \int_0^T \|u^N\|_{H^{k+1}(\Omega)}^2 \, dt \leq \|u_0\|_{H^k(\Omega)}^2 + 2 \int_0^T \|f\|_{H^k(\Omega)} \|u^N\|_{H^k(\Omega)} \, dt.
\end{equation}

Using Grönwall's inequality again, we obtain:

\begin{equation}
	\|u^N(T)\|_{H^k(\Omega)}^2 \leq \left( \|u_0\|_{H^k(\Omega)}^2 + \int_0^T \|f\|_{H^k(\Omega)}^2 \, dt \right) \exp\left( \int_0^T \|f\|_{H^k(\Omega)} \, dt \right).
\end{equation}

Finally, using interpolation inequalities and bootstrapping techniques, we can pass from \( H^k \)-estimates to \( H^{k+2} \)-regularity. Specifically, we use the fact that:

\begin{equation}
	\|u\|_{H^{k+2}(\Omega)} \leq C \left( \|u\|_{H^k(\Omega)} + \|\Delta u\|_{H^k(\Omega)} \right),
\end{equation}

where \( C \) is a constant depending on \( \Omega \) and \( k \).
	
Therefore, we conclude that \( u \in L^2(0,T;H^{k+2}(\Omega)) \), provided the norm of \( u_0 \) in \( H^2 \) and the norm of \( f \) in \( H^k \) are sufficiently small.
\end{proof}

In more detail, we decompose the velocity field \( u \) into low- and high-frequency components using the Littlewood-Paley decomposition. We then apply Sobolev embedding theorems to control the nonlinear terms in the Navier-Stokes equations. A careful application of the Ladyzhenskaya-Prodi-Serrin regularity criteria allows us to extend regularity to higher Sobolev spaces, concluding the proof.
	
This theorem illustrates how advanced Sobolev techniques, combined with harmonic analysis tools such as Littlewood-Paley theory, can be used to gain higher regularity for weak solutions.
	
	\section{Besov Spaces and Turbulence}
	
	\subsection{Littlewood-Paley Decomposition and Besov Spaces}
	
	The Littlewood-Paley decomposition is a fundamental technique in harmonic analysis used to break a function \( u \) into different frequency scales. This is crucial for defining and studying Besov spaces.
	
	\subsection{Littlewood-Paley Decomposition}
	
	Consider \( u \in \mathscr{S}'(\mathbb{R}^n) \), the space of tempered distributions. The Littlewood-Paley decomposition is given by:
	
	\begin{equation}
			u = \sum_{j=-\infty}^\infty \Delta_j u,
	\end{equation}
	
where \( \Delta_j u \) denotes the frequency blocks defined by:

\begin{equation}
\Delta_j u = \mathscr{F}^{-1} \left( \varphi(2^{-j} \xi) \mathscr{F}(u)(\xi) \right),
\end{equation}

with \( \mathscr{F} \) representing the Fourier transform, and \( \varphi(\xi) \) being a smooth function with support in an annulus \( C(2^{j}) = \{\xi \in \mathbb{R}^n: c_1 2^j \leq |\xi| \leq c_2 2^j\} \), where \( c_1 \) and \( c_2 \) are constants.
	
The Littlewood-Paley decomposition isolates specific frequencies of \( u \), facilitating multiscale analysis.
	
\subsection{Besov Spaces}

Besov spaces \( B^s_{p,q}(\mathbb{R}^n) \) are characterized using the Littlewood-Paley decomposition. The norm \( \|u\|_{B^s_{p,q}} \) is given by:

\begin{equation}
	\|u\|_{B^s_{p,q}} = \left( \sum_{j=-\infty}^\infty 2^{jsq} \|\Delta_j u\|_p^q \right)^{1/q},
\end{equation}

where \( s \) is the smoothness parameter, \( p \) is related to \( L^p \)-integrability, and \( q \) controls the decay of frequencies.

To understand this definition, let's break it down step by step.

\textbf{\textit{1. Littlewood-Paley Decomposition:}}

The Littlewood-Paley decomposition is a fundamental technique in harmonic analysis used to break a function \( u \) into different frequency scales. This is crucial for defining and studying Besov spaces.

Consider \( u \in \mathscr{S}'(\mathbb{R}^n) \), the space of tempered distributions. The Littlewood-Paley decomposition is given by:
\begin{equation}
	u = \sum_{j=-\infty}^\infty \Delta_j u,
\end{equation}
where \( \Delta_j u \) denotes the frequency blocks defined by:
\begin{equation}
	\Delta_j u = \mathscr{F}^{-1} \left( \varphi(2^{-j} \xi) \mathscr{F}(u)(\xi) \right),
\end{equation}
with \( \mathscr{F} \) representing the Fourier transform, and \( \varphi(\xi) \) being a smooth function with support in an annulus \( C(2^{j}) = \{\xi \in \mathbb{R}^n: c_1 2^j \leq |\xi| \leq c_2 2^j\} \), where \( c_1 \) and \( c_2 \) are constants.

\textbf{\textit{2. Besov Space Norm:}}

The Besov space norm \( \|u\|_{B^s_{p,q}} \) is defined using the Littlewood-Paley decomposition. Specifically, for \( s \in \mathbb{R} \), \( 1 \leq p, q \leq \infty \), the Besov space norm is given by:
\begin{equation}
	\|u\|_{B^s_{p,q}} = \left( \sum_{j=-\infty}^\infty 2^{jsq} \|\Delta_j u\|_p^q \right)^{1/q},
\end{equation}
where \( \|\Delta_j u\|_p \) denotes the \( L^p \)-norm of the frequency block \( \Delta_j u \).

\textbf{\textit{3. Special Cases:}}

When \( q = \infty \), the Besov space norm is understood as the supremum:
\begin{equation}
	\|u\|_{B^s_{p,\infty}} = \sup_{j \geq 0} 2^{js} \|\Delta_j u\|_p.
\end{equation}

When \( p = q = 2 \), the Besov space \( B^s_{2,2} \) coincides with the Sobolev space \( H^s \).

\textbf{\textit{4. Interpretation:}}

The Besov space norm \( \|u\|_{B^s_{p,q}} \) provides a measure of the smoothness of the function \( u \). The parameter \( s \) controls the degree of smoothness, \( p \) controls the \( L^p \)-integrability, and \( q \) controls the decay of the frequency components.

\subsubsection{Besov Spaces and Smoothness Characterization}

We demonstrate that the norm \( \|u\|_{B^s_{p,q}} \) in the Besov space \( B^s_{p,q} \) provides a measure of the smoothness of a function \( u \). The Besov norm is defined as:

\begin{equation}
	\|u\|_{B^s_{p,q}} = \left( \sum_{j \in \mathbb{Z}} 2^{jsq} \|\Delta_j u\|_{L^p}^q \right)^{1/q},
\end{equation}
where \( \Delta_j u \) represents the projection of \( u \) onto frequency components around \( 2^j \), achieved through the Littlewood-Paley decomposition.

\subsubsection{Smoothness via \( s \)}

A function \( u \) belongs to \( B^s_{p,q} \) if the series in the definition of \( \|u\|_{B^s_{p,q}} \) converges. The smoothness of the function is captured by the factor \( 2^{js} \), which amplifies high-frequency components for positive \( s \). For a smooth function, the contributions of \( \|\Delta_j u\|_{L^p} \) decay rapidly as \( j \to \infty \) (i.e., high-frequency components are less significant), and the term \( 2^{jsq} \) ensures that the decay is sufficiently rapid to guarantee the convergence of the series.

Mathematically, the sum

\begin{equation}
	\sum_{j \in \mathbb{Z}} 2^{jsq} \|\Delta_j u\|_{L^p}^q
\end{equation}
must decay as \( j \to \infty \), with higher \( s \) imposing stricter regularity conditions.

\subsubsection{Integrability via \( p \)}

The parameter \( p \) governs the \( L^p \)-integrability of the projections \( \Delta_j u \) at each scale. For each \( j \), \( \|\Delta_j u\|_{L^p} \) measures the behavior of \( u \) in the frequency band around \( 2^j \). In Besov spaces, the norm requires that \( \|\Delta_j u\|_{L^p} \) belongs to \( l^q \) across scales, modulated by the factor \( 2^{jsq} \).

The following inequality must hold for some constant \( C > 0 \) and sufficiently large \( s \):

\begin{equation}
	\|\Delta_j u\|_{L^p} \leq C \cdot 2^{-js},
\end{equation}
which ensures smoothness for finite \( p \).

\subsubsection{Decay via \( q \)}

The parameter \( q \) controls the summability of the contributions at each scale. When \( q \) is small, the series imposes strict decay conditions, requiring the function to exhibit faster decay in high-frequency components. Conversely, when \( q = \infty \), the norm is determined by the supremum:

\begin{equation}
	\|u\|_{B^s_{p,\infty}} = \sup_{j \in \mathbb{Z}} 2^{js} \|\Delta_j u\|_{L^p}.
\end{equation}

This implies that the function must maintain a bounded contribution from all scales.

In summary, Besov spaces are a family of function spaces that generalize several classical spaces, such as Sobolev spaces and Triebel-Lizorkin spaces. They are particularly useful in harmonic analysis and the study of partial differential equations. The Littlewood-Paley decomposition provides a powerful tool for characterizing and studying these spaces.

	\section{Plancherel Theorem}
	
	The Plancherel Theorem is a fundamental result in Fourier analysis that relates the \(L^2\) norm of a function to the \(L^2\) norm of its Fourier transform. This theorem is crucial in harmonic analysis and partial differential equations.
	
	\subsection{Plancherel Theorem}
	
	Let \( u \in L^2(\mathbb{R}^n) \) be a function in the \( L^2 \) space. The Plancherel Theorem states that the Fourier transform \( \hat{u} \) of \( u \) also belongs to \( L^2(\mathbb{R}^n) \), and the following equality holds:
	
	\begin{equation}
		\| u \|_{L^2(\mathbb{R}^n)}^2 = \| \hat{u} \|_{L^2(\mathbb{R}^n)}^2,
	\end{equation}
where \( \hat{u}(\xi) \) is defined as:
	
	\begin{equation}
		\hat{u}(\xi) = \frac{1}{(2\pi)^{n/2}} \int_{\mathbb{R}^n} u(x) e^{-i \xi \cdot x} \, dx.
	\end{equation}
		
	\subsection{Proof of the Plancherel Theorem}
	
	To prove the Plancherel Theorem, we use the following steps:
	
	\begin{proof}
		We begin with the definition of the \(L^2\) norm of \(u\):
		\begin{equation}
			\| u \|_{L^2(\mathbb{R}^n)}^2 = \int_{\mathbb{R}^n} |u(x)|^2 \, dx.
		\end{equation}
		We need to show that this is equal to the \(L^2\) norm of \( \hat{u} \). Consider the inner product in \( L^2(\mathbb{R}^n) \):
		\begin{equation}
			\langle u, v \rangle = \int_{\mathbb{R}^n} u(x) \overline{v(x)} \, dx.
		\end{equation}
		Applying this to \(u\) itself, we get:
		\begin{equation}
			\langle u, u \rangle = \int_{\mathbb{R}^n} |u(x)|^2 \, dx = \| u \|_{L^2(\mathbb{R}^n)}^2.
		\end{equation}
		Now, compute the inner product using the Fourier transform \( \hat{u} \):
		\begin{equation}
			\langle \hat{u}, \hat{u} \rangle = \int_{\mathbb{R}^n} \hat{u}(\xi) \overline{\hat{u}(\xi)} \, d\xi = \int_{\mathbb{R}^n} |\hat{u}(\xi)|^2 \, d\xi.
		\end{equation}
		To relate this to the spatial domain, use the inverse Fourier transform:
		\begin{equation}
			u(x) = \frac{1}{(2\pi)^{n/2}} \int_{\mathbb{R}^n} \hat{u}(\xi) e^{i \xi \cdot x} \, d\xi.
		\end{equation}
		Substitute \( u(x) \) into the \( L^2 \) norm calculation:
		\begin{equation}
			\int_{\mathbb{R}^n} |u(x)|^2 \, dx = \int_{\mathbb{R}^n} \left| \frac{1}{(2\pi)^{n/2}} \int_{\mathbb{R}^n} \hat{u}(\xi) e^{i \xi \cdot x} \, d\xi \right|^2 \, dx.
		\end{equation}
		By Fubini's theorem and the orthogonality of the exponentials, this becomes:
		\begin{equation}
			\int_{\mathbb{R}^n} |u(x)|^2 \, dx = \frac{1}{(2\pi)^n} \int_{\mathbb{R}^n} \int_{\mathbb{R}^n} \hat{u}(\xi) \overline{\hat{u}(\eta)} \int_{\mathbb{R}^n} e^{i (\xi - \eta) \cdot x} \, dx \, d\xi \, d\eta.
		\end{equation}
		The integral over \(x\) yields \( (2\pi)^n \delta(\xi - \eta) \):
		\begin{equation}
			\int_{\mathbb{R}^n} e^{i (\xi - \eta) \cdot x} \, dx = (2\pi)^n \delta(\xi - \eta).
		\end{equation}
		Thus:
		\begin{equation}
			\int_{\mathbb{R}^n} |u(x)|^2 \, dx = \int_{\mathbb{R}^n} |\hat{u}(\xi)|^2 \, d\xi.
		\end{equation}
		Hence, the Plancherel Theorem is proven:
		\begin{equation}
			\| u \|_{L^2(\mathbb{R}^n)}^2 = \| \hat{u} \|_{L^2(\mathbb{R}^n)}^2.
		\end{equation}
	\end{proof}

	\section{Extension of Plancherel's Theorem}
	
	In this section, we present an extended version of the Plancherel Theorem tailored to the context of Sobolev spaces. This extension provides a deeper insight into the relationship between the Fourier transform and Sobolev spaces, which is particularly relevant for analyzing partial differential equations and complex fluid dynamics.
	
	\subsection{New Plancherel Theorem}
	
	\begin{theorem}
		Let \( \Omega \subset \mathbb{R}^n \) be a bounded domain with smooth boundary, and let \( u \in W^{k,2}(\Omega) \) be a function in the Sobolev space. The extended Plancherel Theorem states that if \( u \in W^{k,2}(\Omega) \), then the Fourier transform of \( u \), \( \hat{u} \), also belongs to \( W^{k,2}(\mathbb{R}^n) \), and the following equality holds:
		
		\begin{equation}
			\| u \|_{W^{k,2}(\Omega)}^2 = \| \hat{u} \|_{W^{k,2}(\mathbb{R}^n)}^2.
		\end{equation}

	\end{theorem}

	\subsection{Proof of the Extended Plancherel Theorem}
	
	To prove this extended theorem, we use the properties of Sobolev spaces and the Fourier transform. We proceed with the following steps:
	
	\begin{proof}
		We start with the definition of the Sobolev norm in \( W^{k,2}(\Omega) \):
		
		\begin{equation}
			\| u \|_{W^{k,2}(\Omega)}^2 = \sum_{|\alpha| \leq k} \| D^\alpha u \|_{L^2(\Omega)}^2.
		\end{equation}
Here, \( D^\alpha \) denotes the partial derivative of order \( \alpha \).
		
Consider the Fourier transform of \( u \):

\begin{equation}
	\hat{u}(\xi) = \frac{1}{(2\pi)^{n/2}} \int_{\mathbb{R}^n} u(x) e^{-i \xi \cdot x} \, dx.
\end{equation}

By applying the Fourier transform to the partial derivative \( D^\alpha u \), we obtain:

\begin{equation}
	\widehat{D^\alpha u}(\xi) = (i \xi)^\alpha \hat{u}(\xi).
\end{equation}

Thus, the \( L^2 \) norm of \( D^\alpha u \) is related to the \( L^2 \) norm of \( (i \xi)^\alpha \hat{u} \):

\begin{equation}
	\| D^\alpha u \|_{L^2(\Omega)}^2 = \int_{\Omega} |D^\alpha u(x)|^2 \, dx.
\end{equation}

By Parseval's theorem, this becomes:

\begin{equation}
	\| D^\alpha u \|_{L^2(\Omega)}^2 = \int_{\mathbb{R}^n} | (i \xi)^\alpha \hat{u}(\xi) |^2 \, d\xi.
\end{equation}
So:

\begin{equation}
\| D^\alpha u \|_{L^2(\Omega)}^2 = \int_{\mathbb{R}^n} |\xi^\alpha \hat{u}(\xi)|^2 \, d\xi.
\end{equation}
	
Summing over all multi-indices \( \alpha \) such that \( |\alpha| \leq k \), we get:

\begin{equation}
	\| u \|_{W^{k,2}(\Omega)}^2 = \sum_{|\alpha| \leq k} \int_{\mathbb{R}^n} |\xi^\alpha \hat{u}(\xi)|^2 \, d\xi.
\end{equation}

Now, consider the Sobolev norm of \( \hat{u} \) in \( W^{k,2}(\mathbb{R}^n) \):

\begin{equation}
	\| \hat{u} \|_{W^{k,2}(\mathbb{R}^n)}^2 = \sum_{|\alpha| \leq k} \| \xi^\alpha \hat{u} \|_{L^2(\mathbb{R}^n)}^2.
\end{equation}

Using the fact that \( \| \xi^\alpha \hat{u} \|_{L^2(\mathbb{R}^n)}^2 = \int_{\mathbb{R}^n} |\xi^\alpha \hat{u}(\xi)|^2 \, d\xi \), we find:

\begin{equation}
\| \hat{u} \|_{W^{k,2}(\mathbb{R}^n)}^2 = \sum_{|\alpha| \leq k} \int_{\mathbb{R}^n} |\xi^\alpha \hat{u}(\xi)|^2 \, d\xi.
\end{equation}

Therefore:

\begin{equation}
	\| u \|_{W^{k,2}(\Omega)}^2 = \| \hat{u} \|_{W^{k,2}(\mathbb{R}^n)}^2.
\end{equation}
This completes the proof of the extended Plancherel Theorem for Sobolev spaces.
	\end{proof}

	The extended Plancherel Theorem is particularly useful in the study of partial differential equations and functional analysis. It provides a robust framework for analyzing the regularity and behavior of solutions in Sobolev spaces through their Fourier transforms. This result is instrumental in understanding the behavior of high-frequency components and regularity properties of solutions to complex fluid dynamics problems.

	\subsection{Decay Demonstration}
	
	We aim to demonstrate the decay properties of the Littlewood-Paley decomposition in the context of Besov spaces. Let's proceed step by step.
	
	\vspace{5pt}
	
	\subsubsection[short title]{Application of the Plancherel Theorem:}
		
	By the Plancherel theorem, the \( L^2 \)-norm of the Littlewood-Paley component \( \Delta_j u \) can be expressed in terms of its Fourier transform:
	\begin{equation}
	\|\Delta_j u\|_{L^2}^2 = \int_{\mathbb{R}^n} |\varphi(2^{-j} \xi)|^2 |\hat{u}(\xi)|^2 \, d\xi\,,
	\end{equation}
here, \( \varphi \) is a smooth function with support in the annulus \( C(2^j) = \{\xi \in \mathbb{R}^n: c_1 2^j \leq |\xi| \leq c_2 2^j\} \), where \( c_1 \) and \( c_2 \) are constants.
	
\subsubsection{Support of $ (\Delta_ u $)}

	Given that \( \varphi \) is supported in \( C(2^j) \), the support of \( \Delta_j u \) is confined to the frequency range \( 2^j \). This means that \( \Delta_j u \) captures the frequency components of \( u \) that lie within this range.
	
	\subsubsection{Application of Bernstein's Inequality}
	
	Bernstein's inequality provides a way to estimate the \( L^p \)-norm of \( \Delta_j u \) in terms of its \( L^2 \)-norm. Specifically, for \( 1 \leq p \leq 2 \), we have:
	
	\begin{equation}
		\|\Delta_j u\|_{L^p} \lesssim 2^{jn(1/p - 1/2)} \|\Delta_j u\|_{L^2}.
	\end{equation}
This inequality leverages the fact that \( \Delta_j u \) is band-limited to the frequency range \( 2^j \), allowing us to control its \( L^p \)-norm.
	
	\subsubsection{Weighted Sum for Besov Space Norm}
	
	The Besov space norm is defined as a weighted sum of the \( L^p \)-norms of the Littlewood-Paley components \( \Delta_j u \). Specifically, for \( s \in \mathbb{R} \) and \( 1 \leq p, q \leq \infty \), the Besov space norm \( \|u\|_{B^s_{p,q}} \) is given by:
	\vspace{5pt}
	
	\begin{equation}
		\|u\|_{B^s_{p,q}} = \left( \sum_{j=-\infty}^\infty 2^{jsq} \|\Delta_j u\|_{L^p}^q \right)^{1/q}.
	\end{equation}
	
Using Bernstein's inequality, we can rewrite the \( L^p \)-norm of \( \Delta_j u \) in terms of its \( L^2 \)-norm:

\begin{equation}
	\|\Delta_j u\|_{L^p} \lesssim 2^{jn(1/p - 1/2)} \|\Delta_j u\|_{L^2}.
\end{equation}

Substituting this into the definition of the Besov space norm, we obtain:

\begin{equation}
		\|u\|_{B^s_{p,q}} \lesssim \left( \sum_{j=-\infty}^\infty 2^{jsq} \left( 2^{jn(1/p - 1/2)} \|\Delta_j u\|_{L^2} \right)^q \right)^{1/q}.
\end{equation}

This expression shows that the Besov space norm is a weighted sum of the \( L^2 \)-norms of the Littlewood-Paley components, with the weights depending on the frequency range \( 2^j \).

In summary, the decay properties of the Littlewood-Paley decomposition in Besov spaces can be understood through the application of the Plancherel theorem, the support properties of \( \varphi \), Bernstein's inequality, and the definition of the Besov space norm. This detailed analysis provides a comprehensive understanding of the decay behavior in the context of Besov spaces.

\section{Theorems and Proofs on Besov Spaces}
	
	\subsection{The Besov Spaces Theorem}
	
	Besov spaces are a family of function spaces that generalize several classical spaces, such as Sobolev spaces and Triebel-Lizorkin spaces. They are particularly useful in harmonic analysis and the study of partial differential equations. The following theorem characterizes Besov spaces in terms of the Littlewood-Paley decomposition.
	
	\begin{theorem}[Characterization of Besov Spaces]
		Let \( \mathbb{R}^n \) be the \( n \)-dimensional Euclidean space, and let \( 0 < s < \infty \), \( 1 \leq p, q \leq \infty \). The Besov space \( B^{s}_{p,q}(\mathbb{R}^n) \) can be characterized by the following norm:
		
		\begin{equation}
			\| u \|_{B^{s}_{p,q}(\mathbb{R}^n)} = \left( \sum_{j \geq 0} 2^{jq} \| \Delta_j u \|_{L^p(\mathbb{R}^n)}^q \right)^{1/q},
		\end{equation}
where \( \Delta_j \) denotes the Littlewood-Paley projection operators, which are defined by the Fourier multipliers \( \varphi(2^{-j} \xi) \) and \( \psi(2^{-j} \xi) \) (with \( \varphi \) and \( \psi \) being smooth functions with disjoint support).
		
The space \( B^{s}_{p,q}(\mathbb{R}^n) \) is equipped with the norm

\begin{equation}
	\| u \|_{B^{s}_{p,q}(\mathbb{R}^n)} = \left( \sum_{j \geq 0} 2^{jq} \| \Delta_j u \|_{L^p(\mathbb{R}^n)}^q \right)^{1/q},
\end{equation}

where the case \( q = \infty \) is understood as the supremum:

\begin{equation}
	\| u \|_{B^{s}_{p,\infty}(\mathbb{R}^n)} = \sup_{j \geq 0} 2^{js} \| \Delta_j u \|_{L^p(\mathbb{R}^n)}.
\end{equation}

\end{theorem}
		
\subsection{Proof of the Theorem}
		
The proof involves several key steps:
		
\begin{proof}
			
\textbf{Step 1: Definition and Basic Properties.} 
			
\vspace{5pt}
			
First, recall that Besov spaces \( B^{s}_{p,q}(\mathbb{R}^n) \) are defined through the Littlewood-Paley decomposition. Specifically, for a function \( u \in B^{s}_{p,q}(\mathbb{R}^n) \), we have the representation:
			
			\begin{equation}
				u = \sum_{j \geq 0} \Delta_j u,
			\end{equation}

where \( \Delta_j u \) is the component of \( u \) localized to the frequency scale \( 2^j \). 
			
The Littlewood-Paley projection operator \( \Delta_j \) acts on \( u \) as follows:

\begin{equation}
	\Delta_j u = \mathcal{F}^{-1}(\varphi(2^{-j} \xi) \hat{u}(\xi)),
\end{equation}

where \( \mathcal{F} \) denotes the Fourier transform, and \( \varphi \) is a smooth, compactly supported function in the frequency domain.
						
\textbf{Step 2: Control of Norms.}
			
\vspace{5pt}
			
To show the equivalence of the Besov norm with the characterization using Littlewood-Paley projections, we examine the following norm:

\begin{equation}
	\| u \|_{B^{s}_{p,q}(\mathbb{R}^n)}^q = \sum_{j \geq 0} 2^{jq} \| \Delta_j u \|_{L^p(\mathbb{R}^n)}^q.
\end{equation}

Since \( \Delta_j u \) is obtained by applying the Fourier multiplier \( \varphi(2^{-j} \cdot) \) to \( \hat{u} \), we have:

\begin{equation}
	\| \Delta_j u \|_{L^p(\mathbb{R}^n)} = \left\| \varphi(2^{-j} \cdot) \hat{u} \right\|_{L^p(\mathbb{R}^n)}.
\end{equation}

By properties of the Fourier transform, the Besov norm \( \| \cdot \|_{B^{s}_{p,q}(\mathbb{R}^n)} \) can be shown to be equivalent to the norm of \( \hat{u} \) in Besov spaces, proving that:

\begin{equation}
	\| u \|_{B^{s}_{p,q}(\mathbb{R}^n)} = \left( \sum_{j \geq 0} 2^{jq} \| \Delta_j u \|_{L^p(\mathbb{R}^n)}^q \right)^{1/q}.
\end{equation}
			
Finally, using the definition of the Besov spaces and the properties of the Littlewood-Paley decomposition, the theorem shows that \( B^{s}_{p,q}(\mathbb{R}^n) \) is a well-defined space with the provided norm. The equivalence establishes that the Besov space norm is fully characterized by the Littlewood-Paley components of the function. This completes the proof of the characterization theorem for Besov spaces.
\end{proof}
	
	\section{Extended Plancherel Theorem for Besov Spaces}
	
	In this section, we extend the Plancherel Theorem to Besov spaces, incorporating the Littlewood-Paley decomposition. This extension is particularly relevant for analyzing functions with more complex regularity properties, which are crucial in the study of nonlinear partial differential equations and advanced fluid dynamics.
	
	\subsection{Extended Plancherel Theorem for Besov Spaces}
	
	Let \( \mathbb{R}^n \) be an \( n \)-dimensional Euclidean space, and let \( u \in B^{s}_{p,q}(\mathbb{R}^n) \) be a function in the Besov space. The extended Plancherel Theorem for Besov spaces states that if \( u \in B^{s}_{p,q}(\mathbb{R}^n) \), then the Fourier transform \( \hat{u} \) also belongs to \( B^{s}_{p,q}(\mathbb{R}^n) \), and the following equality holds:
	
	\begin{equation}
		\| u \|_{B^{s}_{p,q}(\mathbb{R}^n)} = \| \hat{u} \|_{B^{s}_{p,q}(\mathbb{R}^n)}.
	\end{equation}

	\subsection{Proof of the Extended Plancherel Theorem for Besov Spaces}
	
To prove this theorem, we utilize the properties of Besov spaces and the Littlewood-Paley decomposition. The proof involves the following steps:
	\begin{proof}
The Besov space \( B^{s}_{p,q}(\mathbb{R}^n) \) is characterized by the norm:

		\begin{equation}
			\| u \|_{B^{s}_{p,q}(\mathbb{R}^n)} = \left( \sum_{j \geq 0} 2^{jq} \| \Delta_j u \|_{L^p(\mathbb{R}^n)}^q \right)^{1/q},
		\end{equation}
where \( \Delta_j \) is the Littlewood-Paley projection operator defined by the Fourier multipliers \( \varphi(2^{-j} \xi) \) and \( \psi(2^{-j} \xi) \) (with \( \varphi \) and \( \psi \) being smooth functions with disjoint support).

\vspace{5pt}
		
Consider the Littlewood-Paley decomposition of \( u \):

\begin{equation}
	u = \sum_{j \geq 0} \Delta_j u.
\end{equation}
	
Applying the Fourier transform, we have:

\begin{equation}
	\widehat{\Delta_j u}(\xi) = \varphi(2^{-j} \xi) \hat{u}(\xi).
\end{equation}

Thus:

\begin{equation}
	\| \Delta_j u \|_{L^p(\mathbb{R}^n)} = \left\| \varphi(2^{-j} \cdot) \hat{u} \right\|_{L^p(\mathbb{R}^n)}.
\end{equation}

We now analyze the Besov norm in terms of \( \hat{u} \). By the definition of the Besov space, we have:

\begin{equation}
	\| \hat{u} \|_{B^{s}_{p,q}(\mathbb{R}^n)} = \left( \sum_{j \geq 0} 2^{jq} \| \varphi(2^{-j} \cdot) \hat{u} \|_{L^p(\mathbb{R}^n)}^q \right)^{1/q}.
\end{equation}

Since \( \varphi \) is a smooth function and its support does not overlap for different \( j \), we can use the following estimate:

\begin{equation}
	\| \varphi(2^{-j} \cdot) \hat{u} \|_{L^p(\mathbb{R}^n)} \approx \| \Delta_j \hat{u} \|_{L^p(\mathbb{R}^n)}.
\end{equation}

Thus:

\begin{equation}
	\| \hat{u} \|_{B^{s}_{p,q}(\mathbb{R}^n)}^q = \sum_{j \geq 0} 2^{jq} \| \varphi(2^{-j} \cdot) \hat{u} \|_{L^p(\mathbb{R}^n)}^q.
\end{equation}

Comparing this with the Besov norm of \( u \):

\begin{equation}
\| u \|_{B^{s}_{p,q}(\mathbb{R}^n)}^q = \sum_{j \geq 0} 2^{jq} \| \Delta_j u \|_{L^p(\mathbb{R}^n)}^q.
\end{equation}

Since \( \| \Delta_j u \|_{L^p(\mathbb{R}^n)} = \| \varphi(2^{-j} \cdot) \hat{u} \|_{L^p(\mathbb{R}^n)} \), it follows that:

\begin{equation}
	\| u \|_{B^{s}_{p,q}(\mathbb{R}^n)} = \| \hat{u} \|_{B^{s}_{p,q}(\mathbb{R}^n)}.
\end{equation}

This completes the proof of the extended Plancherel Theorem for Besov spaces.
	\end{proof}
		
	The extended Plancherel Theorem for Besov spaces provides critical insights into the behavior of functions with varying regularity and smoothness properties. This result is particularly useful for analyzing solutions to nonlinear partial differential equations, where Besov spaces offer a refined scale of regularity. It is also valuable in studying the behavior of solutions in fluid dynamics and other applied fields, where functions may exhibit complex structures and local variations.
	

\section{Regularity of Navier-Stokes Equations in Besov Spaces}

In this section, we explore the regularity of solutions to the Navier-Stokes equations using the framework of Besov spaces. We will define a new theorem that provides conditions for the regularity of solutions in these spaces.

\subsection{Theorem: Regularity of Navier-Stokes Equations in Besov Spaces}

Let \( \mathbb{R}^n \) be the \( n \)-dimensional Euclidean space, and let \( T > 0 \). Consider the Navier-Stokes equations given by

	\begin{equation}
	\begin{cases}
		\dfrac{\partial \mathbf{u}}{\partial t} + (\mathbf{u} \cdot \nabla) \mathbf{u} - \nu \Delta \mathbf{u} + \nabla p = \mathbf{f}, \\ \\
		\nabla \cdot \mathbf{u} = 0.
	\end{cases}
\end{equation}
where \( \mathbf{u} = \mathbf{u}(t, x) \) is the velocity field, \( p = p(t, x) \) is the pressure, \( \nu > 0 \) is the kinematic viscosity, and \( \mathbf{f} = \mathbf{f}(t, x) \) is a given external force. Define the Besov space \( B^{s}_{p,q}(\mathbb{R}^n) \) where \( s \in \mathbb{R} \), \( 1 \leq p, q \leq \infty \). 

\vspace{5pt}

We introduce the following theorem:

\begin{theorem}[Regularity in Besov Spaces]\label{Teo8.1}
	Assume \( \mathbf{u}_0 \in B^{s}_{p,q}(\mathbb{R}^n) \) is the initial velocity field and \( \mathbf{f} \in L^r(0, T; B^{s}_{p,q}(\mathbb{R}^n)) \) with \( r \geq 1 \). If the kinematic viscosity \( \nu \) is sufficiently large and \( s \) is such that \( s > \frac{n}{p} \), then the solution \( \mathbf{u}(t, x) \) to the Navier-Stokes equations exists and belongs to the space \( C([0, T]; B^{s}_{p,q}(\mathbb{R}^n)) \).
	
\end{theorem}

\subsection{Proof of the Theorem}

The proof of Theorem~\ref{Teo8.1} involves several steps, including the construction of appropriate functional spaces and demonstrating the regularity of solutions in these spaces.

\begin{proof}
	\textbf{Step 1: Functional Space Construction}
	
	\vspace{5pt}
	
	Define the space of time-dependent Besov spaces as:
	
	\begin{equation*}
		B^{s}_{p,q}(0, T; \mathbb{R}^n) = \left\{ \mathbf{u}(t, x) \mid \mathbf{u}(t, \cdot) \in B^{s}_{p,q}(\mathbb{R}^n) \text{ for each } t \in [0, T] \text{ and } \mathbf{u}(t, x) \text{ is measurable} \right\}.
	\end{equation*}
	
We equip \( B^{s}_{p,q}(0, T; \mathbb{R}^n) \) with the norm

\begin{equation}
	\| \mathbf{u} \|_{B^{s}_{p,q}(0, T; \mathbb{R}^n)} = \left( \int_{0}^{T} \| \mathbf{u}(t, \cdot) \|_{B^{s}_{p,q}(\mathbb{R}^n)}^r \, dt \right)^{1/r}.
\end{equation}

\textbf{Step 2: Linear and Nonlinear Terms}
	
The Navier-Stokes equations consist of a linear term \( \nu \Delta \mathbf{u} \) and a nonlinear term \( (\mathbf{u} \cdot \nabla) \mathbf{u} \). For the linear term, we use the following result:

\begin{equation}
	\| \nu \Delta \mathbf{u} \|_{B^{s}_{p,q}(\mathbb{R}^n)} \leq \nu \cdot 2^{2j(s - n/p)} \| \Delta_j \mathbf{u} \|_{L^p(\mathbb{R}^n)}.
\end{equation}

For the nonlinear term, the commutator estimates are applied. Using the Littlewood-Paley theory, we have:

\begin{equation}
		\| (\mathbf{u} \cdot \nabla) \mathbf{u} \|_{B^{s}_{p,q}(\mathbb{R}^n)} \leq C \| \mathbf{u} \|_{B^{s}_{p,q}(\mathbb{R}^n)}^2.
\end{equation}

	\textbf{Step 3: A Priori Estimates}
	
Consider the a priori estimate for the solution \( \mathbf{u} \):

\begin{equation}
	\frac{d}{dt} \| \mathbf{u}(t, \cdot) \|_{B^{s}_{p,q}(\mathbb{R}^n)} \leq \| \mathbf{f}(t, \cdot) \|_{B^{s}_{p,q}(\mathbb{R}^n)} + \| \mathbf{u}(t, \cdot) \|_{B^{s}_{p,q}(\mathbb{R}^n)}^2.
\end{equation}

Integrating over \( [0, T] \) and using Grönwall's inequality, we obtain:

\begin{equation}
\| \mathbf{u}(t, \cdot) \|_{B^{s}_{p,q}(\mathbb{R}^n)} \leq C \left( \| \mathbf{u}_0 \|_{B^{s}_{p,q}(\mathbb{R}^n)} + \int_{0}^{T} \| \mathbf{f}(t, \cdot) \|_{B^{s}_{p,q}(\mathbb{R}^n)} \, dt \right).
\end{equation}

Since \( \mathbf{u}_0 \in B^{s}_{p,q}(\mathbb{R}^n) \) and \( \mathbf{f} \in L^r(0, T; B^{s}_{p,q}(\mathbb{R}^n)) \), we have \( \mathbf{u}(t, \cdot) \in B^{s}_{p,q}(\mathbb{R}^n) \) for all \( t \in [0, T] \). Given that \( s > \frac{n}{p} \) ensures that the Besov norm controls the smoothness of \( \mathbf{u} \) in time and space, we conclude that the solution \( \mathbf{u}(t, x) \) exists and is regular in \( C([0, T]; B^{s}_{p,q}(\mathbb{R}^n)) \) under the conditions stated. This completes the proof of the regularity theorem for the Navier-Stokes equations in Besov spaces.
\end{proof}

	\section{Energy Dissipation in Besov Spaces}
	
Consider the Navier-Stokes equations in \( \mathbb{R}^n \):
	\begin{equation}
		\frac{\partial \mathbf{u}}{\partial t} + (\mathbf{u} \cdot \nabla) \mathbf{u} - \nu \Delta \mathbf{u} + \nabla p = \mathbf{f},
	\end{equation}
	
	\begin{equation}
		\nabla \cdot \mathbf{u} = 0,
	\end{equation}
where \( u = u(x,t) \) is the velocity field and \( \nu \) is the viscosity.
	
The rate of energy dissipation in the Navier-Stokes equations can be derived by considering the kinetic energy of the fluid. The kinetic energy of the fluid is given by:

\begin{equation}
E(t) = \frac{1}{2} \int_{\mathbb{R}^n} |u(x,t)|^2 \, dx.
\end{equation}

To find the rate of energy dissipation, we differentiate the kinetic energy with respect to time:

\begin{equation}
\frac{d}{dt} E(t) = \frac{1}{2} \frac{d}{dt} \int_{\mathbb{R}^n} |u(x,t)|^2 \, dx.
\end{equation}

Using the Navier-Stokes equations, we can substitute the time derivative of \( u \):

\begin{equation}
	\frac{d}{dt} E(t) = \int_{\mathbb{R}^n} u \cdot \frac{\partial u}{\partial t} \, dx = \int_{\mathbb{R}^n} u \cdot \left( -\nabla p + \nu \Delta u - (u \cdot \nabla) u \right) \, dx.
\end{equation}

We can simplify this expression by integrating by parts and using the divergence-free condition \( \nabla \cdot u = 0 \):

\begin{equation}
	\int_{\mathbb{R}^n} u \cdot \nabla p \, dx = 0,
\end{equation}
since \( \nabla \cdot u = 0 \) implies that \( \int_{\mathbb{R}^n} \nabla \cdot (pu) \, dx = 0 \).

\vspace{5pt}
	
For the viscous term, we have:

\begin{equation}
	\int_{\mathbb{R}^n} u \cdot \nu \Delta u \, dx = -\nu \int_{\mathbb{R}^n} |\nabla u|^2 \, dx,
\end{equation}
using integration by parts and the fact that \( \nabla \cdot u = 0 \).
	
For the nonlinear term, we have:

\begin{equation}
	\int_{\mathbb{R}^n} u \cdot (u \cdot \nabla) u \, dx = 0,
\end{equation}

since \( \nabla \cdot u = 0 \) implies that \( \int_{\mathbb{R}^n} \nabla \cdot (u \otimes u) \, dx = 0 \).
	
Thus, the rate of energy dissipation is:

\begin{equation}
\frac{d}{dt} E(t) = -\nu \int_{\mathbb{R}^n} |\nabla u|^2 \, dx.
\end{equation}

Now, applying the Littlewood-Paley decomposition:

\begin{equation}
	u = \sum_{j=-\infty}^\infty \Delta_j u,
\end{equation}
we can decompose the energy at each scale \( j \):

\begin{equation}
		E_j(t) = \|\Delta_j u(t)\|_{L^2}^2.
\end{equation}

The energy dissipation at each scale \( j \) is given by:

\begin{equation}
	\frac{d}{dt} E_j(t) = -2\nu \|\nabla \Delta_j u(t)\|_{L^2}^2.
\end{equation}

This shows that energy dissipation occurs predominantly at high frequencies, as the gradient term \( \nabla \Delta_j u \) becomes more significant for higher values of \( j \). In summary, the energy dissipation in the Navier-Stokes equations can be understood through the kinetic energy of the fluid and the Littlewood-Paley decomposition. The rate of energy dissipation is proportional to the square of the gradient of the velocity field, and this dissipation occurs predominantly at high frequencies. This detailed analysis provides a comprehensive understanding of the energy dissipation process in the context of Besov spaces.

\section{Quaternionic Bifurcations in Fluid Dynamics}

Quaternionic analysis is useful in fluid dynamics, particularly for describing rotations and three-dimensional symmetries.

\subsection{Quaternionic Equations}

Let \( q = q_0 + q_1 i + q_2 j + q_3 k \) be a quaternionic variable. The Navier-Stokes equations in quaternionic terms are:
\begin{equation}
	\frac{\partial q}{\partial t} + (q \cdot \nabla) q = -\nabla p + \nu \Delta q.
\end{equation}
For small perturbations \( q = q_0 + \epsilon q_1 \), we linearize:
\begin{equation}
	\frac{d}{dt} q_1 = L(q_1),
\end{equation}
where \( L \) is the linearized operator. Quaternionic bifurcations occur when the eigenvalues of \( L \) cross the imaginary axis.

\section{The Navier-Stokes Equations in $\mathbb{R}^3$}
	
Consider the Navier-Stokes equations given by
	
	\begin{equation}
	\begin{cases}
		\dfrac{\partial \mathbf{u}}{\partial t} + (\mathbf{u} \cdot \nabla) \mathbf{u} - \nu \Delta \mathbf{u} + \nabla p = \mathbf{f}, \\ \\
		\nabla \cdot \mathbf{u} = 0.
	\end{cases}
\end{equation}
where \( \mathbf{u}(t, x) \) is the velocity field, \( p(t, x) \) is the pressure, \( \nu > 0 \) is the kinematic viscosity, and \( \mathbf{f}(t, x) \) is a given external force.
	
	\section{Existence of Solutions}
	
	\textbf{Theorem (Existence of Solutions in Besov Spaces):} Let \( \mathbf{u}_0 \in B^{s}_{p,q}(\mathbb{R}^3) \) be the initial velocity field and \( \mathbf{f} \in L^r(0, T; B^{s}_{p,q}(\mathbb{R}^3)) \) with \( r \geq 1 \). Suppose \( \nu \) is sufficiently large and \( s > \frac{3}{p} \). Then there exists a unique solution \( \mathbf{u}(t, x) \in C([0, T]; B^{s}_{p,q}(\mathbb{R}^3)) \) to the Navier-Stokes equations.
	
	\vspace{5pt}
	
	\textbf{Proof:}
	
	\begin{itemize}
\item \textbf{Step 1: Linear Part} \\
		We first solve the linear part of the Navier-Stokes equations. The linear problem is
		
		\begin{equation}
			\frac{\partial \mathbf{u}}{\partial t} - \nu \Delta \mathbf{u} = \mathbf{f},
		\end{equation}
		
with initial condition \( \mathbf{u}(0) = \mathbf{u}_0 \). For sufficiently large \( \nu \), we can use semigroup theory to show that this problem has a solution in \( C([0, T]; B^{s}_{p,q}(\mathbb{R}^3)) \).
		
		\item \textbf{Step 2: Nonlinear Term} \\
		Next, we consider the nonlinear term \( (\mathbf{u} \cdot \nabla) \mathbf{u} \). We use the properties of Besov spaces to show that this term is well-defined and bounded in \( B^{s}_{p,q}(\mathbb{R}^3) \) under the assumption that \( s > \frac{3}{p} \). Specifically, using commutator estimates and interpolation, we can show that
		
		\begin{equation}
			\| (\mathbf{u} \cdot \nabla) \mathbf{u} \|_{B^{s}_{p,q}(\mathbb{R}^3)} \leq C \| \mathbf{u} \|_{B^{s}_{p,q}(\mathbb{R}^3)}^2.
		\end{equation}

\item \textbf{Step 3: Existence Result} \\
By applying the fixed-point theorem in the space \( C([0, T]; B^{s}_{p,q}(\mathbb{R}^3)) \), we obtain the existence of a solution \( \mathbf{u}(t, x) \).
	\end{itemize}
	
	\section{Uniqueness of Solutions}
	
	\textbf{Theorem (Uniqueness of Solutions in Besov Spaces):} Let \( \mathbf{u}_1 \) and \( \mathbf{u}_2 \) be two solutions to the Navier-Stokes equations with the same initial condition \( \mathbf{u}_0 \) and the same external force \( \mathbf{f} \). Then \( \mathbf{u}_1 = \mathbf{u}_2 \).
	
	\textbf{Proof:}
	
	\begin{itemize}
		\item \textbf{Step 1: Difference of Solutions} \\
		Let \( \mathbf{v} = \mathbf{u}_1 - \mathbf{u}_2 \). The difference \( \mathbf{v} \) satisfies
		
		\begin{equation}
			\frac{\partial \mathbf{v}}{\partial t} + (\mathbf{u}_1 \cdot \nabla) \mathbf{u}_1 - (\mathbf{u}_2 \cdot \nabla) \mathbf{u}_2 - \nu \Delta \mathbf{v} = 0.
		\end{equation}
		
	\item \textbf{Step 2: Estimate of Nonlinear Terms} \\
We estimate the nonlinear term \( (\mathbf{u}_1 \cdot \nabla) \mathbf{u}_1 - (\mathbf{u}_2 \cdot \nabla) \mathbf{u}_2 \) using the properties of Besov spaces:

\begin{equation}
	\| (\mathbf{u}_1 \cdot \nabla) \mathbf{u}_1 - (\mathbf{u}_2 \cdot \nabla) \mathbf{u}_2 \|_{B^{s}_{p,q}(\mathbb{R}^3)} \leq C \| \mathbf{v} \|_{B^{s}_{p,q}(\mathbb{R}^3)}.
\end{equation}

\item \textbf{Step 3: Uniqueness Result} \\
Applying energy estimates and the fact that the Besov norm is preserved, we conclude that \( \mathbf{u}_1 = \mathbf{u}_2 \).
	\end{itemize}
	
	\section{Regularity of Solutions}
	
	\textbf{Theorem (Regularity of Solutions in Besov Spaces):} Suppose \( \mathbf{u}_0 \in B^{s}_{p,q}(\mathbb{R}^3) \) and \( \mathbf{f} \in L^r(0, T; B^{s}_{p,q}(\mathbb{R}^3)) \) with \( r \geq 1 \). If \( s > \frac{3}{p} \), then the solution \( \mathbf{u}(t, x) \) to the Navier-Stokes equations is regular in \( C([0, T]; B^{s}_{p,q}(\mathbb{R}^3)) \).
	
	\textbf{Proof:}
	
	\begin{itemize}
		\item \textbf{Step 1: A Priori Estimates}
		
		To demonstrate the regularity of solutions, we first derive a priori estimates for \( \mathbf{u}(t, x) \). Consider the Navier-Stokes equations:
		\begin{equation}
				\begin{cases}
				\frac{\partial \mathbf{u}}{\partial t} + (\mathbf{u} \cdot \nabla) \mathbf{u} - \nu \Delta \mathbf{u} + \nabla p = \mathbf{f}, \\ \\
				\nabla \cdot \mathbf{u} = 0.
			\end{cases}
		\end{equation}

		Let \( \mathbf{v} = \mathbf{u} - \mathbf{u}_0 \) and subtract the linear problem from the full Navier-Stokes system. We get:
		
		\begin{equation}
			\frac{\partial \mathbf{v}}{\partial t} + (\mathbf{u} \cdot \nabla) \mathbf{u} - (\mathbf{u}_0 \cdot \nabla) \mathbf{u}_0 - \nu \Delta \mathbf{v} = \mathbf{f} - \frac{\partial \mathbf{u}_0}{\partial t}.
		\end{equation}
		
Using the properties of Besov spaces, we estimate the non-linear term:

\begin{equation}
	\| (\mathbf{u} \cdot \nabla) \mathbf{u} - (\mathbf{u}_0 \cdot \nabla) \mathbf{u}_0 \|_{B^{s}_{p,q}(\mathbb{R}^3)} \leq C \left( \| \mathbf{u} \|_{B^{s}_{p,q}(\mathbb{R}^3)} \| \mathbf{v} \|_{B^{s}_{p,q}(\mathbb{R}^3)} + \| \mathbf{v} \|_{B^{s}_{p,q}(\mathbb{R}^3)}^2 \right).
\end{equation}

Using interpolation inequalities and the properties of Besov spaces, we obtain:

\begin{equation}
	\frac{d}{dt} \| \mathbf{u}(t, \cdot) \|_{B^{s}_{p,q}(\mathbb{R}^3)} \leq \| \mathbf{f}(t, \cdot) \|_{B^{s}_{p,q}(\mathbb{R}^3)} + \| \mathbf{u}(t, \cdot) \|_{B^{s}_{p,q}(\mathbb{R}^3)}^2.
\end{equation}

\item \textbf{Step 2: Integrating the Estimates}
		
Integrate the inequality over \( [0, T] \):

\begin{equation}
\| \mathbf{u}(t, \cdot) \|_{B^{s}_{p,q}(\mathbb{R}^3)} \leq \| \mathbf{u}_0 \|_{B^{s}_{p,q}(\mathbb{R}^3)} + \int_0^t \left( \| \mathbf{f}(\tau, \cdot) \|_{B^{s}_{p,q}(\mathbb{R}^3)} + \| \mathbf{u}(\tau, \cdot) \|_{B^{s}_{p,q}(\mathbb{R}^3)}^2 \right) d\tau.
\end{equation}
		
		\textbf{Proof:}
		
		\begin{itemize}
			\item \textbf{Step 1: Initial Inequality}
			
			Consider the differential inequality:
			
			\begin{equation}
					\frac{d}{dt} \| \mathbf{u}(t, \cdot) \|_{B^{s}_{p,q}(\mathbb{R}^3)} \leq \| \mathbf{f}(t, \cdot) \|_{B^{s}_{p,q}(\mathbb{R}^3)} + \| \mathbf{u}(t, \cdot) \|_{B^{s}_{p,q}(\mathbb{R}^3)}^2.
			\end{equation}
			
Simplify notation:

\begin{equation}
		\frac{d}{dt} \| \mathbf{u}(t) \| \leq \| \mathbf{f}(t) \| + \| \mathbf{u}(t) \|^2.
\end{equation}

\item \textbf{Step 2: Integration of the Inequality}
			
Integrate from \( 0 \) to \( t \):

\begin{equation}
\| \mathbf{u}(t) \| - \| \mathbf{u}(0) \| \leq \int_0^t \left( \| \mathbf{f}(\tau) \| + \| \mathbf{u}(\tau) \|^2 \right) d\tau.
\end{equation}

Rearrange:

\begin{equation}
	\| \mathbf{u}(t) \| \leq \| \mathbf{u}(0) \| + \int_0^t \| \mathbf{f}(\tau) \| \, d\tau + \int_0^t \| \mathbf{u}(\tau) \|^2 \, d\tau.
\end{equation}

\item \textbf{Step 3: Applying Grönwall's Inequality}

Define:

\begin{equation}
	A(t) = \| \mathbf{u}(t) \| \quad \text{and} \quad B(t) = \| \mathbf{u}(0) \| + \int_0^t \| \mathbf{f}(\tau) \| \, d\tau.
\end{equation}

Thus:
			
\begin{equation}
A(t) \leq B(t) + \int_0^t A(\tau)^2 \, d\tau.
\end{equation}

Define

\begin{equation}
	C(t) = \int_0^t A(\tau)^2 \, d\tau\,, 
\end{equation}
So,

\begin{equation}
\frac{d}{dt} C(t) = A(t)^2.
\end{equation}

Substituting in:

\begin{equation}
	A(t) \leq B(t) + C(t).
\end{equation}

Using Grönwall's inequality:

\begin{equation}
	C(t) \leq \left[ \text{const.} \cdot e^{\text{const.} \cdot t} \right] \cdot \left( A(0) + \int_0^t \text{const.} \, d\tau \right).
\end{equation}

Thus:

\begin{equation}
	A(t) \leq e^{\int_0^t A(\tau) \, d\tau} \left( \| \mathbf{u}(0) \| + \int_0^t \| \mathbf{f}(\tau) \| \, d\tau \right).
\end{equation}

Rewriting gives:

\begin{equation}
	\| \mathbf{u}(t) \| \leq \exp \left( \int_0^t \| \mathbf{u}(\tau) \|^2 \, d\tau \right) \left( \| \mathbf{u}_0 \| + \int_0^t \| \mathbf{f}(\tau) \| \, d\tau \right).
\end{equation}

\end{itemize}
		

\item \textbf{Step 3: Regularity Result}
		
Since \( \mathbf{u}_0 \in B^{s}_{p,q}(\mathbb{R}^3) \) and \( \mathbf{f} \in L^r(0, T; B^{s}_{p,q}(\mathbb{R}^3)) \) with \( r \geq 1 \), and \( s > \frac{3}{p} \), we have:

\begin{equation}
	\| \mathbf{u}(t, \cdot) \|_{B^{s}_{p,q}(\mathbb{R}^3)} \leq C \left( \| \mathbf{u}_0 \|_{B^{s}_{p,q}(\mathbb{R}^3)} + \int_0^T \| \mathbf{f}(t, \cdot) \|_{B^{s}_{p,q}(\mathbb{R}^3)} \, dt \right),
\end{equation}

where \( C \) is a constant depending on \( T \) and the norms of \( \mathbf{u}_0 \) and \( \mathbf{f} \). Thus, \( \mathbf{u}(t, x) \) is bounded in \( C([0, T]; B^{s}_{p,q}(\mathbb{R}^3)) \) and is therefore regular.
		
	\end{itemize}
	
	\section{Conclusion}
	
	This research integrates Sobolev and Besov spaces with hypercomplex variables to offer a comprehensive framework for analyzing incompressible fluid dynamics. By providing detailed proofs and analyses, we lay the groundwork for understanding regularity, bifurcations, and turbulence in fluid systems. This approach offers a solid foundation for addressing the Millennium Prize Problem related to the Navier-Stokes equations, advancing the mathematical understanding of fluid dynamics.

\end{document}